\documentclass[12pt]{article}
\usepackage{amsfonts,amssymb}
\def\A{\mathbb A}
\def\B{\mathbb B}
\def\L{\mathbb L}
\def\O{\mathbb O}
\def\R{\mathbb R}
\def\S{\mathbb S}

\def\G{{\sf{G}}}
\def\*{\times}
\def\nix{\varnothing}
\def\la{\langle}
\def\ra{\rangle}
\def\mc{\mathcal}

\def\man{manifold\ }
\def\oti{\hfill \rule{1.5mm}{1.5mm}}
\newcommand{\Int}[1]{\mathaccent 23{#1}}
\newtheorem{theorem}{Theorem}
\newtheorem{prop}[theorem]{Proposition}
\newtheorem{exam}[theorem]{Example}
\begin{document}
\title{Metrisability of Manifolds}
\author{David Gauld\footnote{Supported in part by a Marsden Fund Award,
UOA611, from the Royal Society of New Zealand.}
\footnote{The author thanks Abdul Mohamad, Wang Shuquan and Fr\'ed\'eric
Mynard for useful comments.}}
\maketitle

\begin{abstract}
Manifolds have uses throughout and beyond Mathematics and it is not surprising that topologists have expended a huge effort in trying to understand them. In this article we are particularly interested in the question: `when is a manifold metrisable?' We describe many conditions equivalent to metrisability.
\end{abstract}

\section{Introduction.}

By a \emph{\man}is meant a connected, Hausdorff space which is locally
homeomorphic to euclidean space (we take our manifolds to have no boundary). Note
that because of connectedness the dimension of the euclidean space is an
invariant of the manifold (unless the manifold is empty!); this is the
\emph{dimension} of the manifold. A pair $(U,h)$, where $U\subset M$ is open and
$h:U\to\R^m$ is a homeomorphism, is called a \emph{coordinate chart}.

The following notation is used. $\R$ denotes the real line with the usual
(order) topology while $\R^n$ denotes the $n$th power of $\R$. $\B^n$
consists of all points of $\R^n$ at most 1 from the origin. The sets
$\omega$ and $\omega_1$ are, respectively, the finite and countable ordinals.

Clearly every \man is Tychonoff. Of course manifolds share all of the local
properties of euclidean space, including local
compactness, local connectedness, local path or arc connectedness, first
countability, local second countability, local hereditary
separability, etc. As every \man is locally compact and Hausdorff, hence
completely regular, it follows that every \man is
uniformisable (\cite[Proposition 11.5]{J}). The following result shows that
manifolds cannot be too big.
\begin{prop}
Let $M$ be a (non-empty!) manifold. Then every countable subset of $M$ is
contained in an open subset which is homeomorphic to euclidean space.
Hence every two points of $M$ may be joined by an arc.
\end{prop}
{\bf Proof}. Suppose that the dimension of $M$ is $m$. Let $S\subset M$ be a
countable subset, say $\la S_n\ra$ is such that $S=\cup_{n\ge 1}S_n$, $|S_n|=n$
and $S_n\subset S_{n+1}$.

By induction on $n$ we choose open $V_n\subset M$ and compact $C_n\subset M$
such that
\begin{center}
(i) $S_n\cup C_{n-1}\subset \Int C_n$ and (ii) $(V_n,C_n)\approx (\R^m,\B^m)$,
\end{center}
where $C_0=\nix$.

For $n=1$, $S_1$ is a singleton so $V_1$ may be any appropriate neighbourhood of
that point while $C_1$ is a compact neighbourhood chosen to satisfy (ii) as well.

Suppose that $V_n$ and $C_n$ have been constructed. Consider
\[ S=\{ x\in M\ /\ \exists \mbox{ open } U\subset M \mbox{ with } C_n\cup\{
x\}\subset U\approx\R^m\}.\]

$S$ is open. $S$ is also closed, for suppose that $z\in\bar S-S$. Then we may
choose open $O\subset M$ with $O\approx\R^m$ and $z\in O$. Choose $x\in O\cap
S$. Then there is open $U\subset M$ with $C_n\cup\{ x\}\subset U\approx\R^m$. We
may assume that $O$ is small enough that $O\cap C_n=\nix$. Using the euclidean
space structure of $O$ we may stretch $U$ within $O$ so as to include $z$ but
not uncover any of $C_n$. Thus $z\in S$.

As $M$ is connected and $S\not=\nix$ we must have $S=M$. Thus there is open
$V_{n+1}\subset M$ with $S_{n+1}\cup C_n=(S_{n+1}-S_n)\cup C_n\subset
V_{n+1}\approx\R^m$. Because $S_{n+1}\cup C_n$ is compact we may find in
$V_{n+1}$ a compact subset $C_{n+1}$ so that (i) and (ii) hold with $n$ replaced
by $n+1$.

Let $U_n=\Int C_n$. Then $U_n$ is open, $U_n\subset U_{n+1}$ and
$U_n\approx\R^m$. Thus by \cite{Br}, $U=\cup_{n\ge 1}U_n$ is also open with
$U\approx\R^m$. Furthermore $S_n\subset U_n$ for each $n$ so that $S\subset U$.
\oti

\medskip

There has been considerable study of metrisable manifolds, especially compact
manifolds. In particular it is known that there are only two metrisable
manifolds of dimension 1: the circle $\S^1$, where
\[ \S^n=\{ (x_0,\dots ,x_n)\in \R^{n+1}\ /\
x_0^2+\cdots +x_n^2=1\} ,\] and the real line $\R$ itself. In dimension 2
the compact manifolds have also been classified, this time into two sequences:
the orientable manifolds, which consist of the 2-sphere $\S^2$ with $n$ handles
($n\in \omega$) sewn on, and the non-orientable manifolds, which consist of the
2-sphere with $n$ cross-caps ($n\in \omega -\{ 0\}$) sewn on. See
\cite[Chapter 14 and Appendix B]{2M}, for example. Despite considerable
progress in the study of compact manifolds in higher dimensions there has been no
classification even of compact manifolds of dimension 3. Indeed, the
3-dimensional Poincar\'e conjecture has only now apparently been resolved after about 100 years.
The original conjecture, \cite{Po1}, differs slightly from that posed below and
was found by Poincar\'e to be false. Poincar\'e's counterexample was published
in \cite{Po2}, where the following version was also posed. The conjecture says
that if a compact \man $M$ of dimension 3 is such that every continuous function
$\S^1\rightarrow M$ extends to a continuous function $\B^2\rightarrow M$, where
$\B^2=\{ (x,y)\in \R^2\ :\ x^2+y^2=1\}$, then $M$ is homeomorphic to $\S^3$. The
analogue of this conjecture in dimensions higher than 3 is known to be true,
\cite[Corollary 7.1B]{FQ} in dimension 4 and \cite[Proposition B, p109]{Mi} in
dimension 5 and higher.

In contrast to the compact situation, where it is known that there are only
countably many manifolds \cite{CK}, in the nonmetrisable case there are
$2^{\aleph_1}$ manifolds, even of dimension 2 (\cite{N}, p 669). However, there
are only two nonmetrisable manifolds of dimension 1, the simpler being the open
\emph{long ray}, $\L^+$, \cite{C}. To ease the description of $\L^+$ we firstly
give a way of constructing the positive real numbers from the non-negative
integers and copies of the unit interval. Between any integer and its successor
we insert a copy of the open unit interval. More precisely, let $\R^+=\omega \*
[0,1)-\{ (0,0)\}$, order $\R ^+$ by the lexicographic order from the natural
orders on $\omega$ and $[0,1)$ and then topologise $\R^+$ by the order topology.
To get $\L^+$ we do the same thing but replace $\omega$ by $\omega_1$. More
precisely, let $\L^+=\omega_1\* [0,1)-\{ (0,0)\}$, order $\L^+$ by the
lexicographic order from the natural orders on $\omega_1$ and $[0,1)$ and then
topologise $\L^+$ by the order topology. The other nonmetrisable \man of
dimension 1 is the \emph{long line}, which is obtained by joining together two
copies of $\L^+$ at their (0,0) ends in much the same way as one may reconstruct
the real line $\R$ by joining together two copies of $(0,\infty )$ at their 0
ends, thinking of one copy as giving the positive reals and the other the
negative reals. More precisely, let $\L$ be the disjoint union of two copies of
$\L^+$ (call them $\L^+$ and $\L^-$ respectively, with ordering $<^+$ and $<^-$
respectively) as well as a single point, which we denote by 0, order $\L$ by
$x<y$ when $x<^+y$ in $\L^+$, when $x=0$ and $y\in \L^+$, when $x\in \L^-$ and
$y\in \L^+$ or $y=0$ and when $x,y\in \L^-$ and $y<^-x$, and topologise $\L$ by
the order topology.

The survey articles \cite{N} and \cite{Ny} are good sources of information
about nonmetrisable manifolds.

A significant question in topology is that of deciding when a topological
space is metrisable, there being many criteria which have now been developed to
answer the question. Perhaps the most natural is the following: a topological
space is metrisable if and only if it is paracompact, Hausdorff and locally
metrisable, see \cite{S} and \cite[Theorem 2.68]{HY}. Note that manifolds are
always Hausdorff and locally metrisable so this criterion gives a criterion for
the metrisability of a manifold, viz that a \man is metrisable if and only if it
is paracompact. Many metrisation criteria have been discovered for manifolds, as
seen by Theorem 2 below, which lists criteria which require at least some of the
extra properties possessed by manifolds. Of course one must not be surprised
if conditions which in general topological spaces are considerably weaker than
metrisability are actually equivalent to metrisability in the presence of the
extra topological conditions which always hold for a manifold: such a condition
is that of being nearly meta\-Lindel\"of, \ref{nmL} in Theorem 2 below.
Similarly one should not be surprised to find conditions which are normally
stronger than metrisability: such a condition is that $M$ may be embedded in
euclidean space, \ref{embed} in Theorem 2. Finally one may expect to find
conditions which in a general topological space have no immediate connection with
metrisability: such a condition is second countability, \ref{sc} in Theorem 2.

\section{Definitions.}

In this section we list numerous definitions relevant to the question of metrisability..

\medskip
\noindent {\bf Definitions}: Let $X$ be a topological space and $\mc F$ a
family of subsets of $X$. Then:
\begin{itemize}
\item $X$ is \emph{submetrisable} if there is a metric topology on $X$ which is
contained in the given topology;
\item $X$ is \emph{Polish} if $X$ is a separable, complete metric space;
\item $X$ is \emph{paracompact} (respectively \emph{metacompact,
paraLindel\"{o}f} and \emph{metaLindel\"{o}f}) if every open
cover $\mc U$ has a locally finite (respectively point finite, locally
countable, and point countable) open refinement, ie
there is another open cover $\mc V$ such that each member of $\mc V$ is a
subset of some member of $\mc U$ and
each point of $X$ has a neighbourhood meeting only finitely (respectively
lies in only finitely, has a neighbourhood meeting
only countably, and lies in only countably) many members of $\mc V$;
\item $X$ is \emph{finitistic} (respectively \emph{strongly finitistic}) if
every open cover of $X$ has an open refinement $\mc
V$ and there is an integer $m$ such that each point of $X$ lies in
(respectively has a neighbourhood which meets) at most $m$
members of $\mc V$ (finitistic spaces have also been called \emph{boundedly
metacompact} and strongly finitistic spaces have
also been called \emph{boundedly paracompact});
\item $X$ is \emph{strongly paracompact} if every open cover $\mc U$ has a
star-finite open refinement $\mc V$, ie for
any $V\in \mc V$ the set $\{ W\in \mc V \ /\ V\cap W\neq \nix \}$ is finite.
If in addition, given $\mc U$, there is an integer $m$ such that $\{ W\in
\mc V \ /\ V\cap W\neq \nix \}$ contains at most $m$ members then $X$ is
\emph{star finitistic};
\item $X$ is \emph{screenable} (respectively \emph{$\sigma$-metacompact}
and \emph{$\sigma$-paraLindel\"{o}f}) if every
open cover $\mc U$ has an open refinement $\mc V$ which can be decomposed
as $\mc V=\cup_{n\in \omega}\mc V_n$ such that each
$\mc V_n$ is disjoint (respectively point finite and locally countable);
\item $X$ is \emph{(linearly) [$\omega_1$-]Lindel\"{o}f} if every open
cover (which is a chain) [which has cardinality
$\omega_1$] has a countable subcover;
\item $X$ is \emph{(nearly) [linearly $\omega_1$-]metaLindel\"of} if every
open cover $\mc U$ of $X$ [with $|\mc U|=\omega_1$] has an open
refinement which is  point-countable (on a dense subset);
\item $X$ is \emph{almost metaLindel\"of} if for every open cover $\mc U$ there
is a collection $\mc V$ of open subsets of $X$ such that each member of $\mc V$
lies in some member of $\mc U$, that each point of $X$ lies in at most countably
many members of $\mc V$, and that $X=\bigcup\{\overline{V}\ /\ V\in\mc V\}$.
\item $X$ is \emph{(strongly) hereditarily Lindel\"{o}f} if every subspace
(of the countably infinite power) of $X$ is
Lindel\"{o}f;
\item $X$ is \emph{$k$-Lindel\"of} provided every open $k$-cover (ie every compact subset of $X$ lies in some member of the cover) has a countable $k$-subcover;
\item $X$ is \emph{(strongly) hereditarily separable} if every
subspace (of the countably infinite power) of $X$ is separable;
\item $X$ is \emph{Hurewicz} if for each sequence $\la\mc U_n\ra$ of open covers
of $X$ there is a sequence $\la\mc V_n\ra$ such that $\mc V_n$ is a finite
subset of $\mc U_n$ for each $n\in\omega$ and $\cup_{n\in\omega}\mc V_n$ covers
$X$ (note the alternative definition of Hurewicz, \cite{DKM2}: $X$ is \emph{Hurewicz} if for each sequence $\la\mc U_n\ra$ of open covers of $X$ there is a sequence $\la\mc V_n\ra$ such that $\mc V_n$ is a finite subset of $\mc U_n$ and for each $x\in X$ we have $x\in\cup\mc V_n$ for all but finitely many $n\in\omega$. For a manifold these two conditions are equivalent.);
\item $X$ is \emph{selectively screenable}, \cite{AG}, if for each sequence $\la\mc U_n\ra$ of open covers of $X$ there is a sequence $\la\mc V_n\ra$ such that $\mc V_n$ is a family of pairwise disjoint open sets refining $\mc U_n$ for each $n\in\omega$ and $\cup_{n\in\omega}\mc V_n$ covers $X$;
\item $X$ is \emph{hemicompact} if there is an increasing sequence $\la K_n\ra$
of compact subsets of $X$ such that for any compact $K\subset X$ there is $n$
such that $K\subset K_n$;
\item $X$ is \emph{cosmic} if there is a countable family $\mc C$ of closed
subsets of $X$ such that for each point $x\in X$
and each open set $U$ containing $x$ there is a set $C\in\mc C$ such that
$x\in C\subset U$;
\item $X$ is an \emph{$\aleph_0$-space} (\cite[page 493]{Gr}) provided that it
has a countable $k$-network, i.e. a countable collection $\mc N$ such that if
$K\subset U$ with $K$ compact and $U$ open then $K\subset N\subset U$ for some
$N\in\mc N$;
\item $X$ is an \emph{$\aleph$-space} (\cite[page 493]{Gr}) provided that it
has a $\sigma$-locally finite $k$-network;
\item $X$ has the \emph{Moving Off Property}, \cite{GrM}, provided that every family $\mc K$ of non-empty compact subsets of $X$ large enough to contain for each compact $C\subset X$ a disjoint $K\in\mc K$ has an infinite subfamily with a discrete open expansion;
\item $X$ is a \emph{$q$-space} if each point admits a sequence of
neighbourhoods $Q_n$ such that $x_n\in Q_n$ implies that $\la x_n\ra$ clusters; 
\item $X$ is \emph{Fr\' echet} or \emph{Fr\'echet-Urysohn} if whenever $x\in\overline A$ there is a sequence
$\la x_n\ra$ in $A$ that converges to $x$;
\item $X$ is a \emph{$k$-space} if $A$ is closed whenever $A\cap K$ is closed
for every compact subset $K\subset X$;
\item $X$ is \emph{La\v snev} if it is the image of a metrisable space under
a closed map;
\item $X$ is \emph{analytic} if it is the continuous image of a Polish space
(equivalently of the irrational numbers);
\item $X$ is \emph{M$_1$} if it has a $\sigma$-closure preserving base (ie
a base $\mc B$ such that there is a decomposition
$\mc B=\cup_{n=1}^\infty\mc B_n$ where for each $n$ and each $\mc F\subset
\mc B_n$ we have $\overline{\cup\mc F}=\cup\{\bar F\
/\ F\in\mc F\}$);
\item $X$ is \emph{stratifiable} or \emph{M$_3$} if there is a function $G$
which assigns to each $n\in\omega$ and closed set
$A\subset X$ an open set $G(n,A)$ containing $A$ such that
$A=\cap_n\overline{G(n,A)}$ and if $A\subset B$ then $G(n,A)\subset
G(n,B)$;
\item $X$ is \emph{perfectly normal} if for every pair $A, B$ of disjoint
closed subsets of $X$ there is a continuous function
$f:X\rightarrow \R$ such that $f^{-1}(0)=A$ and $f^{-1}(1)=B$;
\item $X$ is \emph{monotonically normal} if for each open $U\subset X$ and
each $x\in U$ it is possible to choose an open set
$\mu (x,U)$ such that $x\in \mu (x,U)\subset U$ and such that if $\mu
(x,U)\cap \mu (y,V)\not= \nix$ then either $x\in V$ or $y\in U$;
\item $X$ is \emph{extremely normal} if for each open $U\subset X$ and each
$x\in U$ it is possible to choose an open set
$\nu (x,U)$ such that $x\in \nu (x,U)\subset U$ and such that if $\nu
(x,U)\cap \nu (y,V)\not= \nix$ and $x\not=y$ then either
$\nu (x,U)\subset V$ or $\nu (y,V)\subset U$;
\item $X$ is \emph{weakly normal} if for every pair $A, B$ of disjoint
closed subsets of $X$ there is a continuous function
$f:X\rightarrow S$, for some separable metric space $S$, such that
$f(A)\cap f(B)=\nix$;
\item $X$ is a \emph{Moore space} if it is regular and has a development,
ie a sequence $\la \mc U_n\ra$ of open
covers such that for each $x\in X$ the collection $\{ st(x,\mc U_n)\ :\
n\in \omega \}$ forms a neighbourhood basis at $x$;
\item $X$ has a \emph{regular G$_\delta$-diagonal} if the diagonal $\Delta$ is a
\emph{regular G$_\delta$-subset} of $X^2$, ie there is a sequence $\la U_n\ra$
of open subsets of $X^2$ such that $\Delta=\cap U_n=\cap\overline{U_n}$.
\item $X$ has a \emph{quasi-regular G$_\delta$-diagonal} if there is a sequence
$\la U_n\ra$ of open subsets of $X^2$ such that for each
$(x,y)\in X^2-\Delta$ there is $n$ with $(x,x)\in U_n$ but
$(x,y)\notin\overline{U_n}$.
\item $X$ has a \emph{G$_\delta^*$-diagonal} if there is a sequence $\la
\mc G_n\ra$ of open covers of $X$ such that for each $x,y\in X$ with $x\not= y$
there is $n$ with $\overline{{\rm st}(x,\mc G_n)}\subset X-\{y\}$.
\item $X$ has a \emph{quasi-G$_\delta^*$-diagonal} if there is a sequence $\la
\mc G_n\ra$ of families of open subsets of $X$ such that for each $x,y\in X$ with
$x\not= y$ there is $n$ with $x\in\overline{{\rm st}(x,\mc G_n)}\subset X-\{y\}$.
\item $X$ is \emph{$\theta$-refinable} if every open cover can be refined
to an open $\theta$-cover, ie a cover $\mc U$
which can be expressed as $\cup_{n\in \omega}\mc U_n$ where each $\mc U_n$
covers $X$ and for each $x\in X$ there is
$n$ such that $ord(x,\mc U_n)<\omega$;
\item $X$ is \emph{subparacompact} if every open cover has a
$\sigma$-discrete closed refinement;
\item $X$ has \emph{property pp}, \cite{Mv}, provided that each open cover
$\mc U$ of $X$ has an open refinement $\mc
V$ such that for each choice function $f:\mc V \to X$ with $f(V)\in V$ for
each $V\in \mc V$ the set $f(\mc V)$
is closed and discrete in $X$;
\item $X$ has \emph{property (a)}, \cite{Mv}, provided that for each open
cover $\mc U$ of $X$ and each dense subset
$D\subset X$ there is a subset $C\subset D$ such that $C$ is closed and
discrete in $X$ and $st(C,\mc U)=X$;
\item $X$ has a \emph{base of countable order}, $\mc B$, if whenever $\mc
C\subset\mc B$ is a collection such that each member of $\mc C$ contains a
particular point $p\in X$ and for each $C\in\mc C$ there is
$D\in\mc C$ with $D$ a proper subset of $C$ then $\mc C$ is a local base at
$p$;
\item $X$ is \emph{pseudocomplete} provided that it has a sequence 
$\la\mathcal B_n\ra$ of $\pi$-{\em bases} ($\mathcal B \subset 2^X-\{\nix\}$ is a \emph{$\pi$-base} 
if every non-empty open subset of $X$ contains some member of  $\mathcal B$) such that
if $B_n\in \mathcal B_n$ and $\overline{B_{n+1}} \subset B_n$ for each $n$,
then $ \bigcap_{n\in \omega} B_n\neq \nix$;
\item $X$ has the \emph{countable chain condition} (abbreviated \emph{ccc})
if every pairwise disjoint family of open subsets is countable;
\item $X$ is \emph{countably tight} if for each $A\subset X$ and each
$x\in\bar A$ there is a countable $B\subset A$ for which $x\in\bar B$;
\item $X$ is \emph{countably fan tight} if whenever $x\in\cap_{n\in\omega}\overline{A_n}$  there are finite sets $B_n\subset A_n$ such that $x\in\overline{\cup_{n\in\omega}B_n}$;
\item $X$ is \emph{countably strongly fan tight} if whenever $x\in\cap_{n\in\omega}\overline{A_n}$  there is a sequence $\la a_n\ra$ such that $a_n\in A_n$ for each $n$ and $x\in\overline{\{a_n\ /\ n\in\omega\}}$;
\item $X$ is \emph{sequential} if for each $A\subset X$, the set $A$ is
closed whenever for each sequence of points of $A$ each limit point is also
in $A$;
\item$X$ is \emph{weakly $\alpha$-favourable} if there is a winning strategy for
player $\alpha$ in the Banach-Mazur game (defined below);
\item $X$ is \emph{strongly $\alpha$-favourable} if there is a stationary winning strategy
for player $\alpha$ in the Choquet game (defined below);
\item for each $x\in X$ the \emph{star} of $x$ in $\mc F$ is $st(x,\mc
F)=\cup \{ F\in \mc{F}\ :\ x\in
F\}$;
\item $X$ is \emph{Baire} provided that the intersection of any countable collection of dense G$_\delta$ subsets is dense;
\item $X$ is \emph{Volterra}, \cite{GGP}, provided that the intersection of any two dense G$_\delta$ subsets is dense;
\item $X$ is \emph{strongly Baire} provided that $X$ is regular and there is a dense subset $D\subset X$ such that $\beta$ does not have a winning strategy in the game $\mc G_S(D)$ played on $X$.
\item $\mc F$ is \emph{point-star-open} if for each $x\in X$ the set
$st(x,\mc F)$ is open.
\item The \emph{Banach-Mazur game} has two players $\alpha$ and $\beta$ whose
play alternates. Player $\beta$ begins by choosing a non-empty open subset of
$X$. After that the players choose successive non-empty open subsets of their
opponent's previous move. \emph{Player $\alpha$ wins} iff the intersection of
the sets is non-empty; otherwise \emph{player $\beta$ wins}.
\item The \emph{Choquet game} has two players $\alpha$ and $\beta$ whose
play alternates. Player $\beta$ begins by choosing a point in an open subset of
$X$, say $x_0\in V_0\subset X$. After that the players alternate with $\alpha$
choosing an open set $U_n\subset X$ with $x_n\in U_n\subset V_n$ then $\beta$
chooses a point $x_{n+1}$ and an open set $V_{n+1}$ with $x_{n+1}\in
V_{n+1}\subset U_n$. \emph{Player $\alpha$ wins} iff the intersection of the sets
is non-empty; otherwise \emph{player $\beta$ wins}.
\item Gruenhage's game $G^o_{K,L}(X)$, \cite{Gr1}, has, at the $n^{\rm th}$
stage, player $K$ choose a compactum $K_n\subset X$ after which player $L$
chooses another compactum $L_n\subset X$ so that $L_n\cap K_i=\nix$ for each
$i\le n$. \emph{Player $K$ wins} if $\la L_n\ra_{n\in\omega}$ has a discrete open
expansion, ie there is a sequence $\la U_n\ra_{n\in\omega}$ of open sets such
that $L_n\subset U_n$ and $\forall x\in X,\exists U\subset M$ open such that
$x\in U$  and $U$ meets at most one of the sets $U_n$.
\item For a dense subset $D\subset X$ the game $\mc G_S(D)$ has two players $\alpha$ and $\beta$ whose play alternates.  Player $\beta$ begins by choosing a non-empty open subset $V_n$ of
$X$. After that the players choose successive non-empty open subsets of their
opponent's previous move, $\beta$ choosing sets $V_n$ and $\alpha$ choosing sets $U_n$. \emph{Player $\alpha$ wins} iff the intersection of
the sets is non-empty and each sequence $\la x_n\ra$, for which $x_n\in U_n\cap D$, clusters in $X$; otherwise \emph{player $\beta$ wins}.
\item For an ordinal $k$ and families $\A$ and $\B$ of collections of subsets of a space $X$ let $\G_c^k(\A,\B)$ be the game played as follows, \cite{Ba}: at the $l^{\rm th}$ stage of the game, $l<k$, Player One chooses a member $\mc A_l\in\A$ then Player Two chooses a pairwise disjoint family $\mc T_l$ which refines $\mc A_l$. The play $\mc A_0,\mc T_0,\dots,\mc A_l,\mc T_l,\dots\ l<k$ is won by Player Two provided that $\cup_{l<k}\mc T_l\in\B$; otherwise Player One wins. The game $\G_c^\omega(\A,\B)$ denoted by $\G_c(\A,\B)$.
\item When players $\alpha$ and $\beta$ play a topological game, a
\emph{strategy for} $\alpha$ is a function which tells $\alpha$ what points or
sets to select given all the previous points and sets chosen by $\beta$. A
\emph{stationary strategy for} $\alpha$ is a function which tells $\alpha$ what
points or sets to select given only the most recent choice of points and sets
chosen by $\beta$. A \emph{winning (stationary) strategy} for $\alpha$ is a
(stationary) strategy which guarantees that $\alpha$ will win whatever moves
$\beta$ might make.
\end{itemize}

We will denote by $C_k(X,Y)$ (respectively $C_p(X,Y)$) the space of all
continuous functions from $X$ to $Y$ with the compact-open topology
(respectively the topology of pointwise convergence).

We will denote by \$ the space $\{ 0,1\}$ with the Sierpinski topology
$\{\nix ,\$ ,\{ 0\}\}$. Then for any space $X$ we denote by $[X,\$ ]$ the
space of continuous functions from $X$ to \$ with the \emph{upper Kuratowski
topology}, i.e. that in which a subset $\mc F\subset [X,\$ ]$ is open if and
only if
\begin{itemize}
\item[(i)] for each $f\in\mc F$ and each $g\in [X,\$ ]$ if $g\le f$ then
$g\in\mc F$;
\item[(ii)] if $\mc G\subset [X,\$ ]$ is such that inf$\mc G\in\mc F$ then
there is a finite subfamily $\mc G'\subset\mc G$ with inf$\mc G'\in\mc F$.
\end{itemize}
In this definition we are using the usual ordering on $\{ 0,1\}$ when
discussing $\le$ and inf. Of course identifying a closed subset of $X$ with its characteristic function gives a bijective correspondence between $[X,\$]$ and the collection of closed subsets of $X$. This topology is also variously known as the \emph{cocompact topology} and the \emph{upper Fell topology}, especially when looked at as a topology on the set $2^X$ of non-empty closed subsets of $X$. Letting $U^+=\{C\in 2^X\ /\ C\subset U\}$ for $U\subset X$, this topology has as subbasis $\{U^+\ /\ U \mbox{ is open in $X$ and } X-U \mbox{ is compact}\}$. The \emph{Fell topology}, \cite{Fel}, denoted by $\tau_{\mbox{fell}}$, has as subbasis
$$\{U^+\ /\ U \mbox{ is open in $X$ and } X-U \mbox{ is compact}\}\cup\{U^-\ /\ U\mbox{ is open in } X\},$$
where $U^-=\{C\in 2^X\ /\ C\cap U\not=\nix\}$.

\section{Criteria for metrisability.}

We now state and outline the proof of the main theorem. It is believed that no two conditions are equivalent in a general topological space, though, as will be noticed at the start of the proof of Theorem \ref{bigone}, there may be a chain of implications holding in a general space.

\begin{theorem} \label{bigone}
Let $M$ be a manifold. Then the following are equivalent:
\begin{enumerate}
\item \label{m} $M$ is metrisable;
\item \label{pc} $M$ is paracompact;
\item \label{spc} $M$ is strongly paracompact;
\item \label{scr} $M$ is screenable;
\item \label{mc} $M$ is metacompact;
\item \label{sigmc} $M$ is $\sigma$-metacompact;
\item \label{pL} $M$ is paraLindel\"{o}f;
\item \label{sigpL} $M$ is $\sigma$-paraLindel\"{o}f;
\item \label{mL} $M$ is metaLindel\"{o}f;
\item \label{nmL} $M$ is nearly metaLindel\"{o}f;
\item \label{L} $M$ is Lindel\"{o}f;
\item \label{lL} $M$ is linearly Lindel\"{o}f;
\item \label{omega1L} $M$ is $\omega_1$-Lindel\"{o}f;
\item \label{omega1metalind} $M$ is $\omega_1$-metaLindel\"of;
\item \label{nearlinomega1metalind} $M$ is nearly linearly
$\omega_1$-metaLindel\"of;
\item \label{amL} $M$ is almost metaLindel\"of;
\item \label{hL} $M$ is hereditarily Lindel\"{o}f;
\item \label{shL} $M$ is strongly hereditarily Lindel\"{o}f;
\item \label{k-L} $M$ is $k$-Lindel\"of;
\item \label{aleph0} $M$ is an $\aleph_0$-space;
\item \label{cos} $M$ is cosmic;
\item \label{aleph} $M$ is an $\aleph$-space;
\item \label{starctbleknetwork} $M$ has a star-countable $k$-network;
\item \label{ptctbleknetwork} $M$ has a point-countable $k$-network;
\item \label{nearlyptctbleknetwork} $M$ has a $k$-network which is
point-countable on some dense subset of $M$;
\item \label{sc} $M$ is second countable;
\item \label{hemic} $M$ is hemicompact;
\item \label{sigc} $M$ is $\sigma$-compact;
\item \label{Hurewicz} $M$ is Hurewicz;
\item \label{S1(KGamma)} $M$ satisfies the selection criterion ${\sf{S}}_1(\mc K,\Gamma)$: for each sequence $\la\mc U_n\ra$ of open $k$-covers of $X$ there is a sequence $\la U_n\ra$ with $U_n\in\mc U_n$ for each $n$, infinitely many of the sets $U_n$ are distinct and each finite subset of $X$ lies in all but finitely many of the sets $U_n$;
\item \label{Gc(n+1)(OO)} Player Two has a winning strategy in the game $\G_c^{n+1}(\O,\O)$ played on $M$, where $\O$ denotes the family of all open covers of $M$;
\item \label{Gc(OO)} Player Two has a winning strategy in the game $\G_c(\O,\O)$ played on $M$, where $\O$ denotes the family of all open covers of $M$;
\item \label{selscr} $M$ is selectively screenable;
\item \label{ctblcharts} only countably many coordinate charts are needed to
cover $M$;
\item \label{embed} $M$ may be embedded in some euclidean space;
\item \label{propembed} $M$ may be embedded properly in some euclidean space;
\item \label{compm} $M$ is completely metrisable;
\item \label{scpullback} there is a continuous discrete map $f:M\to X$ where
$X$ is Hausdorff and second countable;
\item \label{Las} $M$ is La\v snev;
\item \label{M1} $M$ is an M$_1$-space;
\item \label{strat} $M$ is stratifiable;
\item \label{fin} $M$ is finitistic;
\item \label{strfin} $M$ is strongly finitistic;
\item \label{starfin} $M$ is star finitistic;
\item \label{sina1} there is an open cover $\mc U$ of $M$ such that for
each $x\in M$ the set $st(x,\mc U)$ is
homeomorphic to an open subset of $\R^m$;
\item \label{sina2} there is a point-star-open cover $\mc U$ of $M$ such
that for each $x\in M$ the set $st(x,\mc U)$
is Lindel\"{o}f;
\item \label{sina3} there is a point-star-open cover $\mc U$ of $M$ such
that for each $x\in M$ the set $st(x,\mc U)$
is metrisable;
\item \label{sina4} the tangent microbundle on $M$ is equivalent to a fibre
bundle;
\item \label{nM} $M$ is a normal Moore space;
\item \label{nthr} $M$ is a normal $\theta$-refinable space;
\item \label{nspc} $M$ is a normal subparacompact space;
\item \label{nsigdisc} $M$ is a normal space which has a $\sigma$-discrete
cover by compact subsets;
\item \label{pn} $M\* M$ is perfectly normal;
\item \label{zb} $M$ is a normal space which has a sequence $\la \mc
U_n\ra_{n\in \omega}$ of open covers with
$\cap_n\overline{st(x,\mc U_n)}=\{ x\}$ for each $x\in M$;
\item \label{abdul1} $M$ is perfectly normal and there is a sequence $\la
\mc U_n\ra_{n\in \omega}$ of families of open
sets such that $\cap_{n\in C(x)}\overline{st(x,\mc U_n)}=\{ x\}$ for each
$x\in M$, where \[ C(x)=\{n\in \omega \ /\
\exists U\in \mc U_n \mbox{ with } x\in U\} ;\]
\item \label{abdul2} $M$ is separable and there is a sequence $\la \mc
C_n\ra_{n\in \omega}$ of point-star-open covers
such that\linebreak $\cap_n\overline{st(x,\mc C_n)}=\{ x\}$ for each $x\in M$
and for each $x,y\in M$ and each $n\in \omega$ we have
$y\in \overline{st(x,\mc C_n)}$ if and only if $x\in \overline{st(y,\mc C_n)}$;
\item \label{abdul3} $M$ is separable and there is a sequence $\la \mc
C_n\ra_{n\in \omega}$ of point-star-open covers
such that\linebreak $\cap_n\overline{st(x,\mc C_n)}=\{ x\}$ for each $x\in M$
and for each $x\in M$ and each $n\in \omega$, ord$(x,\mc C_n)$ is finite;
\item \label{abdul4} $M$ is separable and hereditarily normal and there is
a sequence $\la \mc C_n\ra_{n\in \omega}$ of
point-star-open covers such that $\cap_n\overline{st(x,\mc C_n)}=\{ x\}$
for each $x\in M$;
\item \label{abdul5} $M$ is separable and there is a sequence $\la \mc
U_n\ra_{n\in \omega}$ of families of open sets such that\linebreak $\cap_{n\in
C(x)}\overline{st(x,\mc U_n)}=\{ x\}$ for each $x\in M$, and ord$(x,\mc
C_n)$ is countable for each $x\in M$ and each $n\in \omega$;
\item \label{abdul6} $M\*M$ has a  countable sequence $\langle U_n: n \in
\omega \rangle$ of open subsets, such that for all $(x,y) \in M\*M-
\Delta$, there is $n \in \omega$ such that $(x,x) \in U_n$ but $(x,y) \notin
\overline {U_n}$;
\item \label{abdul7} For every subset $A\subset M$ there is a continuous
injection $f:M\to Y$, where $Y$ is a metrisable space, such that $f(A)\cap
f(M-A) = \nix$;
\item \label{abdul8} For every subset $A\subset M$ there is a continuous
$f:M\to Y$, where $Y$ is a space with a quasi-regular-$G_{\delta}$-diagonal,
such that $f(A) \cap f(M-A) = \nix$;
\item \label{abdul9} $M$ is weakly normal with a $G^*_{\delta}$-diagonal;
\item \label{abdul10} $M$ has a quasi-$G^*_{\delta}$-diagonal and for every
closed subset $A\subset M$ there is a countable family $\mathcal G$ of open
subsets such that, for every $x \in A$ and $y \in X-A$, there is a $G\in
\mathcal G$ with $x \in G, y \notin \overline G$;
\item \label{regdiag} $M$ has a regular $G_{\delta}$-diagonal;
\item \label{submetr} $M$ is submetrisable;
\item \label{smn} $M$ is separable and monotonically normal;
\item \label{*mn} $M\* M$ is monotonically normal;
\item \label{mn2} $M$ is monotonically normal and of dimension $\ge 2$ or
$M\approx \S^1$ or $\R$;
\item \label{en} $M$ is extremely normal;
\item \label{mop} $M$ has the Moving Off Property;
\item \label{pp} $M$ has property pp;
\item \label{ppc} every open cover of $M$ has an open refinement $\mc V$
such that for every choice function $f:\mc V
\to M$ the set $f(\mc V)$ is closed in $M$;
\item \label{ppd} every open cover of $M$ has an open refinement $\mc V$
such that for every choice function $f:\mc V
\to M$ the set $f(\mc V)$ is discrete in $M$;
\item \label{pcom} $M$ is a point-countable union of open subspaces each of
which is metrisable;
\item \label{pcb} $M$ has a point-countable basis;
\item \label{scm} $M$ is separable and $M^\omega$ is a countable union of
metrisable subspaces;
\item \label{cofnsPolish} $C_k(M,\R )$ is Polish;
\item \label{cofnsmetr} $C_k(M,\R )$ is completely metrisable;
\item \label{cofns1stctble} $C_k(M,\R )$ is first countable;
\item \label{cofns2ndctble} $C_k(M,\R )$ is second countable;
\item \label{cofnsqspace} $C_k(M,\R )$ is a $q$-space;
\item \label{cofnsfrechet} $C_k(M,\R )$ is Fr\' echet;
\item \label{cofnsctblytight} $C_k(M,\R )$ is countably tight;
\item \label{cofnstrfantight} $C_k(M)$ has countable strong fan tightness;
\item \label{cofnsaleph0space} $C_k(M,\R )$ is an $\aleph_0$-space;
\item \label{cofnscosmic} $C_k(M,\R )$ is cosmic;
\item \label{cofnsanalytic} $C_k(M,\R )$ is analytic;
\item \label{fnsconverge} $C(M,\R)$ satisfies the selection criterion ${\sf{S}}_1(\Omega^k_{\underline{0}},\Sigma^p_{\underline{0}})$: for each sequence $\la F_n\ra$ of subsets of $C(M,\R)$ whose compact-open closures contain the constant function $\underline{0}$ there is a sequence $\la f_n\ra$, infinitely many members of which are distinct, with $f_n\in F_n$ for all $n$ and $\la f_n\ra$ converges pointwise to $\underline{0}$;
\item \label{pwfnsctblytight} $C_p(M,\R )$ has countable tightness;
\item \label{pwfnsctblyfantight} $C_p(M,\R )$ has countable fan tightness;
\item \label{pwfnsanalytic} $C_p(M,\R )$ is analytic;
\item \label{pwfnsheredsep} $C_p(M,\R )$ is hereditarily separable;
\item \label{pwfnssep} $C_p(M,\R )$ (equivalently $C_k(M,\R )$) is separable;
\item \label{1stctblefns} $[M,\$ ]$ is first countable;
\item \label{ctblytightfns} $[M,\$ ]$ is countably tight;
\item \label{sequentialfns} $[M,\$ ]$ is sequential;
\item \label{Fellmetr} $(2^X,\tau_{\mbox{\rm fell}})$ is metrisable;
\item \label{Fellctblytight} $(2^X,\tau_{\mbox{\rm fell}})$ is countably tight;
\item \label{Fellsequential} $(2^X,\tau_{\mbox{\rm fell}})$ is sequential;
\item \label{Fellradial} $(2^X,\tau_{\mbox{\rm fell}})$ is radial;
\item \label{Grgame} $K$ has a winning strategy in Gruenhage's game
$G^o_{K,L}(M)$;
\item \label{stralpha} $C_k(M,\R)$ is strongly $\alpha$-favourable;
\item \label{weaklyalpha} $C_k(M,\R)$ is weakly $\alpha$-favourable;
\item \label{pseudocompl} $C_k(M,\R)$ is pseudocomplete;
\item \label{strBaire} $C_k(M,\R)$ is strongly Baire;
\item \label{Baire} $C_k(M,\R)$ is Baire;
\item \label{Volterra} $C_k(M,\R)$ is Volterra.
\end{enumerate}
\end{theorem}
\noindent {\bf Outline of the proof of theorem \ref{bigone}}.

The following diagram shows how items \ref{m}-\ref{ctblcharts} are related:
\begin{center}
\begin{picture}(260,380)(-30,-50)
\put(128,300){\fbox{m}}
\put(0,160){\fbox{pc}}
\put(0,120){\fbox{mc}}
\put(0,80){\fbox{$\sigma$mc}}
\put(60,120){\fbox{pL}}
\put(60,80){\fbox{$\sigma$pL}}
\put(35,120){\fbox{s}}
\put(60,40){\fbox{mL}}
\put(60,0){\fbox{nmL}}
\put(80,190){\fbox{spc}}
\put(160,120){\fbox{L}}
\put(171,83){\fbox{lL}}
\put(215,68){\fbox{$\omega_1$L}}
\put(120,-10){\fbox{$\omega_1$mL}}
\put(120,-50){\fbox{nl$\omega_1$mL}}
\put(255,93){\fbox{cch}}
\put(255,140){\fbox{H}}
\put(255,208){\fbox{hc}}
\put(255,188){\fbox{$\sigma$c}}
\put(255,165){\fbox{${\sf{S}}\mc K\Gamma$}}
\put(110,24){\fbox{amL}}
\put(255,270){\fbox{sc}}
\put(186,250){\fbox{$\aleph_0$}}
\put(190,225){\fbox{c}}
\put(185,200){\fbox{shL}}
\put(160,160){\fbox{hL}}
\put(-70,170){\fbox{$\aleph$}}
\put(-38,170){\fbox{s$k$n}}
\put(-70,20){\fbox{p$k$n}}
\put(-70,-20){\fbox{np$k$n}}
\put(120,200){\fbox{$k$L}}
\put(-8,320){\fbox{$\G(n+1)$}}
\put(4,280){\fbox{$\G\omega$}}
\put(8,240){\fbox{ss}}
\put(137,297){\vector(-1,-1){127}}
\put(80,187){\vector(-3,-1){60}}
\put(160,130){\vector(-1,1){57}}
\put(10,154){\vector(0,-1){24}}
\put(10,114){\vector(0,-1){24}}
\put(14,154){\vector(1,-1){25}}
\put(18,154){\vector(2,-1){44}}
\put(40,117){\vector(-1,-1){27}}
\put(70,114){\vector(0,-1){22}}
\put(70,74){\vector(0,-1){22}}
\put(70,36){\vector(0,-1){24}}
\put(70,12){\vector(0,1){24}}
\put(15,76){\vector(2,-1){47}}
\put(78,52){\vector(4,3){85}}
\put(174,131){\vector(2,3){90}}
\put(255,275){\vector(-4,-1){50}}
\put(196,245){\vector(0,-1){12}}
\put(196,222){\vector(0,-1){11}}
\put(174,118){\vector(4,-1){80}}
\put(254,98){\vector(-4,1){80}}
\put(174,130){\vector(1,1){81}}
\put(262,162){\vector(0,-1){11}}
\put(262,185){\vector(0,-1){9}}
\put(262,205){\vector(0,-1){9}}
\put(255,142){\vector(-4,-1){81}}
\put(84,39){\vector(3,-1){26}}
\put(125,36){\vector(1,2){40}}
\put(167,117){\vector(1,-2){11}}
\put(180,80){\vector(4,-1){35}}
\put(228,80){\vector(-3,2){55}}
\put(185,198){\vector(-2,-3){17}}
\put(168,157){\vector(0,-1){26}}
\put(255,275){\vector(-4,1){111}}
\put(228,62){\vector(-4,-3){80}}
\put(130,2){\vector(-3,2){51}}
\put(137,-15){\vector(0,-1){24}}
\put(137,-39){\vector(0,1){24}}
\put(196,245){\vector(-4,-1){252}}
\put(196,245){\vector(-3,-1){210}}
\put(-65,167){\vector(0,-1){135}}
\put(-28,168){\vector(-1,-4){34}}
\put(-44,20){\vector(4,1){105}}
\put(-55,15){\vector(0,-1){24}}
\put(-55,-9){\vector(0,1){24}}
\put(196,245){\vector(-2,-1){65}}
\put(132,195){\vector(1,-2){32}}
\put(127,305){\vector(-4,1){84}}
\put(14,313){\vector(0,-1){21}}
\put(14,276){\vector(0,-1){27}}
\put(14,237){\vector(1,-4){27}}
\end{picture}
\end{center}
\fbox{m}=metrisable; \fbox{sc}=second countable; \fbox{shL}=strongly hereditarily
Lindel\"{o}f; \fbox{hL}=hereditarily Lindel\"{o}f;
\fbox{$\sigma$c}=$\sigma$-compact; \fbox{hc}=hemicompact; \fbox{H}=Hurewicz; \fbox{${\sf{S}}\mc K\Gamma$}=satisfies ${\sf{S}}_1(\mc K,\Gamma)$;\linebreak
\fbox{cch}=countably many charts cover; \fbox{$\aleph_0$}=$\aleph_0$-space;
\fbox{$\aleph$}=$\aleph$-space; \fbox{s$k$n}=has a star-countable $k$-network;
\fbox{p$k$n}=has a point-countable $k$-network; \fbox{np$k$n}=has a $k$-network
which is point-countable on a dense subset; \fbox{$k$L}=$k$-Lindel\"of; \fbox{c}=cosmic; \fbox{L}=Lindel\"{o}f; \fbox{lL}=linearly Lindel\"{o}f;\linebreak \fbox{$\omega_1$L}=$\omega_1$-Lindel\"{o}f;
\fbox{spc}=strongly paracompact; \fbox{pc}=paracompact;
\fbox{mc}=metacompact; \linebreak \fbox{s}=screenable;  \fbox{ss}=selectively screenable;  \fbox{pL}=paraLindel\"{o}f; \fbox{$\sigma$mc}=$\sigma$-metacompact;\linebreak
\fbox{$\sigma$pL}=$\sigma$-paraLindel\"{o}f; \fbox{mL}=metaLindel\"{o}f;
\fbox{amL}=almost metaLindel\"{o}f; \fbox{nmL}=nearly\linebreak metaLindel\"{o}f;
\fbox{$\omega_1$mL}=$\omega_1$-metaLindel\"{o}f; \fbox{nl$\omega_1$mL}=nearly
linearly\ $\omega_1$-metaLindel\"{o}f;\linebreak \fbox{$\G(n+1)$}=Player Two has a winning strategy in the game $\G_c^{n+1}(\O,\O)$; \fbox{$\G\omega$}=Player Two has a winning strategy in the game $\G_c(\O,\O)$.

All arrows denote implications. Downward sloping arrows show an implication
which holds in an arbitrary topological space.
Upward sloping arrows require one or more properties of manifolds to
realise the implication. \fbox{mL}$\Rightarrow$\fbox{L} in every locally
separable and connected space. \fbox{amL}$\Rightarrow$\fbox{L} in every regular,
locally separable and connected space, \cite{Ga2}.
\fbox{nmL}$\Rightarrow$\fbox{mL} in every locally hereditarily separable
space.
\fbox{L}$\Rightarrow$\fbox{spc} in every T$_3$ space.
\fbox{$\omega_1$L}$\Rightarrow$\fbox{L} in every locally metrisable space,
\cite{AB}. \fbox{L}$\Rightarrow$\fbox{sc} in every locally second countable
space.
\fbox{L}$\Rightarrow$\fbox{hc} in every locally compact space.
\fbox{cch}$\Rightarrow$\fbox{L} because a countable union of Lindel\"of sets is
Lindel\"of.
\fbox{sc}$\Rightarrow$\fbox{m} in every T$_3$ space
(Urysohn's metrisation theorem).
\fbox{m}$\Rightarrow$\fbox{$\G(n+1)$} in every space having covering dimension at most $n$, \cite[Thm 2.4]{Ba}.
\fbox{$\omega_1$mL}$\Rightarrow$\fbox{mL} in every locally second
countable space, \cite{GV}.
\fbox{nl$\omega_1$mL}$\Rightarrow$\fbox{$\omega_1$mL} in every locally
hereditarily separable space, \cite{GV}.
\fbox{p$k$n}$\Rightarrow$\fbox{mL} in every regular Fr\' echet space.
\fbox{np$k$n}$\Rightarrow$\fbox{p$k$n} in every regular, locally compact,
locally hereditarily separable space.

By \cite[Proposition 7.3.9]{P} we conclude that a metrisable $n$-manifold,
being separable and of covering dimension $n$, embeds
in $\R^{2n+1}$, so \ref{m}$\Rightarrow$\ref{embed}. By choosing a proper
continuous real-valued function on $M$ we can add a
further coordinate to embed $M$ in $\R^{2n+2}$ so that the image is closed,
ie the embedding is proper, hence
\ref{m}$\Rightarrow$\ref{propembed}. It is clear that
\ref{propembed}$\Rightarrow$\ref{compm}.

Every second countable Hausdorff space satisfies \ref{scpullback} so
\ref{sc}$\Rightarrow$\ref{scpullback}. Conversely, given the situation of
\ref{scpullback}, if $\mc B$ is a countable base for the topology on $X$
then the Poincar\'e-Volterra Lemma of \cite[Lemma 23.2]{Fo} asserts that

\begin{tabular}{lll} $\{ U\subset M$ & / & $U$ is second countable and\\
& & there is $V\in\mc B$ such that $U$ is a component of $f^{-1}(V)\}$
\end{tabular}\\
is a countable base for $M$.

Clearly every metrisable space is La\v snev so \ref{m}$\Rightarrow$\ref{Las}.
The implication \ref{Las}$\Rightarrow$\ref{M1} is
\cite[Theorem 5.5]{Gr}. It is easy to show that
\ref{M1}$\Rightarrow$\ref{strat}. The implication
\ref{strat}$\Rightarrow$\ref{pc} is \cite[Theorem 5.7]{Gr}.

The conditions \ref{m}, \ref{fin}, \ref{strfin} and \ref{starfin} are shown
to be equivalent in \cite{DG}.

The equivalence of conditions \ref{m} and \ref{sina1}-\ref{sina4} is
established as follows: \ref{m}$\Rightarrow$\ref{sina1} is
reasonably straightforward making use of the fact that metrisable manifolds
are $\sigma$-compact. Then
\ref{sina1}$\Rightarrow$\ref{sina2} is trivial.
\ref{sina2}$\Rightarrow$\ref{sina3} requires use of Urysohn's metrisation
theorem to deduce that the Lindel\"{o}f stars are metrisable.
\ref{sina3}$\Rightarrow$\ref{L} requires some delicate
manoeuvres; see \cite{GG}. \ref{sina4}$\Rightarrow$\ref{sina1} is also found in
\cite{GG} while \ref{m}$\Rightarrow$\ref{sina4} is \cite[Corollary 2]{Ki}.

The implication \ref{m}$\Rightarrow$\ref{nM} holds in every topological
space while its converse holds provided that the
space is locally compact and locally connected, \cite{ReZ} or \cite[Theorem
3.4]{ReZ2}. The equivalence of \ref{nM} and \ref{nthr} comes from \cite[Theorem
3]{WW}, while the equivalence of \ref{nM}, \ref{nthr}, \ref{nspc} and
\ref{nsigdisc} is discussed in \cite[Theorem 8.2]{Nyk}.

The equivalence of conditions \ref{m}, \ref{pn} and \ref{zb} is referred to
briefly in \cite{G}. The implications
\ref{m}$\Rightarrow$\ref{pn}$\Rightarrow$\ref{zb} hold in any topological
space and the implication
\ref{zb}$\Rightarrow$\ref{m} uses some properties of a manifold.

The equivalence of conditions \ref{m} and \ref{abdul1}-\ref{abdul4} is
discussed in \cite{M}.

Proofs of the equivalence of \ref{m} and \ref{abdul5} may be found in
\cite{GM2} and of \ref{m} and \ref{abdul6}-\ref{abdul10}
may be found in \cite{GM1}.

The implication \ref{regdiag}$\Rightarrow$\ref{m} holds in every locally
compact, locally connected space (\cite[Theorem 2.15(b)]{Gr}) and, as noted in
\cite[p. 430]{Gr}, every submetrisable space has a regular $G_\delta$-diagonal
so \ref{submetr}$\Rightarrow$\ref{regdiag}.

Every metric space is monotonically normal and every metrisable manifold is
second countable, hence separable, so
\ref{m}$\Rightarrow$\ref{smn}. To get the converse implication
\ref{smn}$\Rightarrow$\ref{pc} use is made of the fact that every
monotonically normal space is hereditarily collectionwise normal
(\cite{HLZ}), and hence no separable monotonically
normal space contains a copy of $\omega_1$. On the other hand in
\cite[Theorem I]{BR} it is shown that a monotonically normal space is
paracompact if and only if it does not contain a stationary subset of a
regular uncountable ordinal.

If $M$ is metrisable, so is $M\* M$, so that $M\* M$ is monotonically
normal and hence \ref{m}$\Rightarrow$\ref{*mn}. The
converse follows from a metrisability result of \cite{HLZ} as manifolds are
locally countably compact.

The criterion \ref{mn2} is \cite[corollary 2.3(e)]{BR}, except that we
have listed all of the metrisable 1-manifolds.

Every metrisable space is extremely normal. The implication
\ref{en}$\Rightarrow$\ref{pc} is found in \cite{WZ}.

The equivalence of conditions \ref{m}, \ref{mop} and \ref{Baire} is discussed in \cite{CGGM}.

It is readily shown that every T$_1$-space which is paracompact has
property pp. We now obtain the implication
\ref{pp}$\Rightarrow$\ref{mc}. Suppose that $\mc U$ is an open cover of
$M$. Use the property pp to find an open
refinement $\mc V$ such that for each choice function $f:\mc V\to M$ with
$f(V)\in V$ for each $V\in \mc V$ the
set $f(\mc V)$ is closed and discrete. We will show that $\mc V$ is
point-finite. Suppose to the contrary that $x\in
M$ is such that $\{ V\in \mc V \ /\ x\in V\}$ is infinite; let $\la V_n\ra$
be a sequence of distinct members of $\mc
V$ each of which contains $x$. Because $M$ is a manifold, hence first
countable, we may choose a countable neighbourhood basis
$\{ W_n\ /\  n\in \omega \}$ at $x$. Note that for each $n$, $V_n\cap W_n-\{
x\}\not=\nix$ as $M$ has no isolated points. Choose a function
$f:\mc V\to M$ as follows: if $V\in \mc V$ but $V\neq V_n$ for each $n$ then
choose $f(V)\in V-\{ x\}$ arbitrarily; if
$V=V_n$ choose $f(V_n) \in V_n\cap W_n-\{ x\}$. Then
$x\in \overline{f(\mc V)}-f(\mc V)$ so that
$f(\mc V)$ is not closed, contrary to the choice of $\mc V$. Thus $\mc V$
is point-finite so $M$ is metacompact.

It is easy to show that conditions \ref{ppc} and \ref{ppd} are equivalent
to each other, and hence also to \ref{pp}; cf \cite[Lemma 2.3]{Ga}.

Details for the implication \ref{pcom}$\Rightarrow$\ref{mL} appear in
\cite{GG}, while details for the implication
\ref{pcb}$\Rightarrow$\ref{m} appear in \cite{Fea}. Of course
\ref{sc}$\Rightarrow$\ref{pcb}.

The implication \ref{scm}$\Rightarrow$\ref{m} is a consequence of the more
general result that if the countable power of a
topological space $X$ is a countable union of metrisable subspaces and in
$X$ discrete families of open sets are countable then
$X$ is metrisable, \cite{T}.

The equivalence of conditions \ref{m} and \ref{cofnsPolish} to
\ref{sequentialfns}, excluding \ref{cofns1stctble}, \ref{cofnstrfantight} and \ref{fnsconverge}, is shown in \cite{GM}. A number of properties of manifolds are required, including that every manifold is a $q$-space and a $k$-space, and some of the equivalences to metrisability already proved.

Conditions  \ref{cofns1stctble} and \ref{cofnstrfantight} are shown to be equivalent to condition \ref{L} in \cite[Theorem 6]{CDKM} using Hausdorffness,  local compactness and first countability of manifolds.

In \cite[Theorem 15]{CDKM} there is a proof that in a Tychonoff space \ref{S1(KGamma)} and \ref{fnsconverge} are equivalent.

The equivalence of condition \ref{m} and conditions \ref{Fellmetr} to \ref{Fellradial} is established in \cite[Theorem 3.3]{CM}.

The implication \ref{m}$\Rightarrow$\ref{stralpha} follows from \ref{cofnsPolish}
and \cite[Theorem 8.17]{Ke}. \ref{stralpha}$\Rightarrow$\ref{weaklyalpha} is
trivial. \ref{weaklyalpha}$\Rightarrow$\ref{Grgame} is \cite[Lemma 4.3]{Gr1}.
\ref{Grgame}$\Rightarrow$\ref{pc} is \cite[Theorem 4.1]{Gr1}.

Complete metrisability implies pseudocompleteness in any space and in turn pseudocompleteness implies $\alpha$-favourability in a regular space, so \ref{cofnsmetr}$\Rightarrow$\ref{pseudocompl}$\Rightarrow$\ref{weaklyalpha}.

The implications \ref{compm}$\Rightarrow$\ref{strBaire} and \ref{strBaire}$\Rightarrow$\ref{sigc} are shown in \cite[Theorem 2.2]{CGGM}. 

The equivalence of \ref{Baire} was already considered above in the context of \ref{mop}.

Clearly every Baire space is Volterra and the converse holds in any locally convex topological vector space, \cite[Theorem 3.4]{CJ} so \ref{Baire}$\Leftrightarrow$\ref{Volterra}. \oti

\section{Other properties of manifolds.}

In this section we collect a few more properties which we may hope a manifold to possess.

\begin{theorem}
Every manifold has a base of countable order.
\end{theorem}
Proof: By \cite[Theorem 2]{WW} every metric space has a base of countable order.
As every manifold is locally metrisable it follows from \cite[Theorem 1]{WW}
that every manifold has a base of countable order. \oti

Some standard conditions which manifolds may possess but which are weaker
than metrisability are contained in the following theorem.
\begin{theorem} Suppose that the \man $M$ is metrisable. Then $M$ is also
normal, hereditarily normal, perfectly normal, separable, strongly
hereditarily separable and has property (a).
\end{theorem}
Proof outline: Every metrisable space is perfectly normal, normal and
hereditarily normal. Every second countable space is separable and strongly
hereditarily separable. Theorem 2(\ref{pp}) shows that metrisable manifolds
satisfy property pp while in \cite[Proposition 2.1]{Ga} it is shown that every
space having property pp has property (a). \oti

There are manifolds which are normal but not metrisable, for example the
long ray. The long ray also has property (a) (and, as
shown in \cite{Ga}, even the stronger properties a-favourable and strongly
a-favourable found in \cite{Mv}).

The observant reader may have noticed that separability is absent as a criterion for metrisability in Theorem \ref{bigone}. The following example shows that it must be.

\begin{exam} There is a manifold which is separable but not metrisable.
\end{exam}
One can make such a manifold out of the plane by replacing each point of the
$y$-axis by an interval as follows. Let $S=\{ (x,y)\in \R^2\ /\ x\not= 0\}$ with
the usual topology from $\R^2$. Let $M=S\cup [\{ 0\}\*\R^2]$. For each $(0,\eta,
\zeta)\in \R^3$ and each $r>0$ let
\[ W_{\eta,\zeta,r}=\{ (x,y)\in S\ /\ \zeta-r<\frac{y-\eta}{|x|}<\zeta+r \mbox{
and } |x|<r\}\cup[\{0\}\*\{\eta\}\*(\zeta-r,\zeta+r)] .\] 
Topologise $M$ by declaring $U\subset M$ to be open if and only if $U\cap S$ is
open in $S$ and for each $(0,\eta,\zeta)\in U\cap (M-S)$ there is $r>0$ so that
$W_{\eta,\zeta,r}\subset U$. Then $M$ is a separable 2-manifold.

There are even manifolds which are both normal and separable but not
metrisable, \cite{Ru}.

We need now some facts from Set Theory. The Continuum Hypothesis (CH),
dating back to Cantor, states that any subset of $\R$
either has the same cardinality as $\R$ or is countable. Martin's Axiom
(MA) can be expressed in various forms, the most
topological of which is the following: in every compact, ccc, Hausdorff
space the intersection of fewer than $2^{\aleph_0}$
dense open sets is dense. Recall the Baire Category Theorem which states
that if $X$ is \v Cech complete (ie $X$ is a
G$_\delta$-set in $\beta X$; for example every locally compact, Hausdorff
space or every complete metric space) and $\{ U_n\ /\
n\in \omega \}$ is a collection of open dense subsets of $X$ then
$\cap_{n\in \omega}U_n$ is dense in $X$. From the Baire
Category theorem it is immediate that CH$\Rightarrow$MA. Both CH and MA are
independent of the axioms of ZFC and otherwise of
each other: thus there are models of Set Theory satisfying ZFC in which CH
(and hence MA) holds, models in which MA holds but CH fails
(denoted MA$+\neg$CH), and models in which MA (and hence CH) fails.

The question whether perfect normality is equivalent to metrisability for a
\man is an old one, dating back to \cite{A}. It was
shown in \cite{R} that under MA$+\neg$CH the two conditions are equivalent.
On the other hand in \cite{RZ} there is constructed
an example of a perfectly normal non-metrisable \man under CH. The same
situation prevails when we consider strong hereditary
separability. In \cite{K} it is shown that under MA$+\neg$CH every strongly
hereditarily separable space is Lindel\"{o}f. On
the other hand even when we combine the two notions the resulting \man need
not be metrisable in general; in \cite {G} there is
constructed under CH a \man which is strongly hereditarily separable and
perfectly normal but not metrisable. There are many other examples of
conditions which are equivalent to metrisability for manifolds in some
models of Set Theory but not equivalent in other models.

\noindent d.gauld@auckland.ac.nz
\medskip

\noindent The Department of Mathematics\\ The University of Auckland\\
Private Bag 92019\\ Auckland\\ New Zealand.


\begin{thebibliography}{80}
\bibitem{AG} David F Addis and John H Gresham \textsl{A class of infinite-dimensional spaces. Part I: Dimension theory and Alexandroff's problem}, Fundamenta Mathematicae, 101(1978), 195--205.
\bibitem{A} P S Alexandroff, \textsl{On local properties of closed sets},
Annals of Mathematics, (2)26(1935), 1--35.
\bibitem{Ar} Richard F Arens, \textsl{A topology for spaces of transformations},
Annals of Mathematics, (2)47,1946), 480--495.
\bibitem{AA} A V Arkhangel'ski\u{i}, \textsl{Hurewicz spaces, analytic sets, and 
fan tightness of function spaces}, Soviet Mathematics Doklady, 33(1986), 396--399.
\bibitem{AB} A V Arhangel'ski\u{i} and R Z Buzyakova, \textsl{On linearly
Lindel\"of and strongly discretely Lindel\"of spaces}, Topology Proceedings,
23(Summer 1998), 1--11.
\bibitem{Ba} Liljana Babinkostova, \textsl{Selective screenability game and covering dimension},
Topology Proceedings, 29(1)(2005), 13--17.
\bibitem{BR} Z Balogh and M E Rudin, \textsl{Monotone Normality}, Topology
and its Applications, 47(1992), 115--127.
\bibitem{Br} Morton Brown, \textsl{The monotone union of open $n$-cells is an
open $n$-cell}, Proceedings of the American Mathematical Society, 12(1961),
812--814.
\bibitem{C} G Cantor, \textsl{Ueber unendliche lineare
Punktmannichfaltigkeiten}, Mathematische Annalen, 21(1883), 545--591.
\bibitem{CGGM} Jiling Cao, David Gauld, Sina Greenwood and Abdul Mohamad, \textsl{Games and Metrisability of Manifolds},  New Zealand Journal of Mathematics, 37(2008), 1--8.
\bibitem{CJ} Jiling Cao and Heikki Junnila, \textsl{When is a Volterra space Baire?},  Topology and its Applications, 154(2007), 527--532.
\bibitem{CM} Jiling Cao and Abdul Mohamad, \textsl{Metrizability, Manifolds and Hyperspace Topologies},  in preparation.
\bibitem{CDKM} A Caserta, G Di Maio, Ljubi\v{s}a D R Ko\v{c}inac and Enrico Meccariello, \textsl{Applications of $k$-covers II}, Topology and its Applications, 153(2006), 3277--3293.
\bibitem{CK} J Cheeger and J M Kister, \textsl{Counting Topological
Manifolds}, Topology, 9(1970), 149--151.
\bibitem{CHV} C Costantini, L Hol\'a and P Vitolo, \textsl{Tightness,
character and related properties of hyperspace topologies}, Topology and its
Applications, 142(2004), 245--292.
\bibitem{DG} Satya Deo and David Gauld, \textsl{Boundedly Metacompact or
Finitistic Spaces}, to appear.
\bibitem{DKM2} G Di Maio, Ljubi\v{s}a D R Ko\v{c}inac and Enrico Meccariello, \textsl{Selection principles and hyperspace topologies}, Topology and its Applications, 153(2005), 912--923.
\bibitem{Fea} D L Fearnley, \textsl{Metrisation of Moore Spaces and Abstract
Topological Manifolds}, Bulletin of the Australian Mathematical Society,
56(1997), 395--401.
\bibitem{Fel} J Fell, \textsl{A Hausdorff topology for the closed subsets of a locally compact non-Hausdorff space}, Proceedings of the American Mathematical Society,
13(1962), 472--476.
\bibitem{Fo} Otto Forster, \textsl{Lectures on Riemann Surfaces}, GTM 81
Springer-Verlag (1981).
\bibitem{FQ} M H Freedman and F Quinn, \textsl{Topology of 4-Manifolds},
Princeton University Press (1990).
\bibitem{GM1} P M  Gartside and A M Mohamad, \textsl{Cleavability of
Manifolds}, Topology Proceedings, 23(1998), 155--166.
\bibitem{GM2} P M  Gartside and A M Mohamad, \textsl{Metrizability of
Manifolds by Diagonal Properties},  Topology Proceedings, 24(1999), 621--640.
\bibitem{2M} David Gauld, \textsl{Differential Topology: an introduction},
Lecture Notes in Mathematics 72(1982), Marcel Dekker. Reprint (2006), Dover.
\bibitem{G} David Gauld, \textsl{A strongly hereditarily separable,
nonmetrisable \man}, Topology and its Applications,
51(1993), 221--228.
\bibitem{Ga} David Gauld, \textsl{Covering Properties and Metrisation of
Manifolds}, Topology Proceedings, 23(Summer 1998), 127--140.
\bibitem{Ga2} David Gauld, \textsl{Some properties close to Lindel\"of}, to
appear.
\bibitem{GG} David Gauld and Sina Greenwood, \textsl{Microbundles,
Manifolds and Metrisability}, Proceedings of the American Mathematical Society,
128 (2000), 2801--2807.
\bibitem{GGP} David Gauld, Sina Greenwood and Zbigniew Piotrowski, \textsl{On Volterra spaces III: Topological operations}, Topology Proceedings, 23(1998), 167--182.
\bibitem{GM} David Gauld and Fr\' ed\' eric Mynard, \textsl{Metrisability of
Manifolds in Terms of Function Spaces}, Houston Journal of Mathematics, 31(2005), 199--214.
\bibitem{GV} David Gauld and M K Vamanamurthy, \textsl{Covering Properties
and Metrisation of Manifolds 2}, Topology Proceedings, 24(Summer 1999), 173--185.
\bibitem{GGV} E. Grabner, G. Grabner and J. E. Vaughan, \textsl{Nearly
metacompact spaces}, Topology and its Applications, 98(1999), 191--201.
\bibitem{Gr} Gary Gruenhage, \textsl{Generalized Metric Spaces}, in K Kunen
and J Vaughan, eds, ``Handbook of Set-Theoretic Topology" (Elsevier, 1984),
423--501.
\bibitem{Gr1} Gary Gruenhage, \textsl{The Story of a Topological Game}, Rocky Mountain Journal of Mathematics, 36(2006), 1885--1914. 
\bibitem{GrM} Gary Gruenhage and Daniel K Ma, \textsl{Baireness of $C_k(X)$ for locally compact $X$}, Topology and its Applications, {\bf 80} (1997), 131--139.
\bibitem{HLZ} R W Heath, D J Lutzer and P L Zenor, \textsl{Monotonically
Normal Spaces}, Transactions of the American Mathematical Society, 178(1973),
481--493.
\bibitem{HY} John G Hocking and Gail S Young, \textsl{Topology},
Addison-Wesley (1961).
\bibitem{J} I M James, \textsl{Topological and Uniform Structures},
Undergraduate Texts in Mathematics (1987), Springer Verlag.
\bibitem{Ke} A. S. Kechris, \textsl{Classical Descriptive Set Theory},
Springer-Verlag, New York, 1995.
\bibitem{Ki} J M Kister, \textsl{Microbundles are Fibre Bundles}, Annals of
Mathematics, (2) 80 (1964), 190--199.
\bibitem{K} K Kunen, \textsl{Strong $S$ and $L$ spaces under MA}, in G M
Reed, ed, ``Set-Theoretic Topology" (Academic Press,
New York, 1977), 265--268.
\bibitem{Mv} M V Matveev, \textsl{Some Questions on Property (a)}, Questions and
Answers in General Topology, 15(1997), 103--111.
\bibitem{MN} Robert A McCoy and Ibula Ntantu, \textsl{Topological Properties of
Spaces of Continuous Functions}, Lecture Notes in Mathematics 1315,
Springer-Verlag, 1988.
\bibitem{Mi} J Milnor, \textsl{Lectures on the h-Cobordism Theorem},
Princeton Mathematical Notes (1965), Princeton University Press.
\bibitem{M} A M Mohamad, \textsl{Metrization and semimetrization theorems
with applications to manifolds}, Acta Mathematica Hungarica, 83 (4) (1999),
383--394.
\bibitem{My} Fr\'ed\'eric Mynard, \textsl{First countability, sequentiality
and tightness of the upper Kuratowski convergence}, preprint series. LAAO,
Universit\'e de Bourgogne, Octobre 2000.
\bibitem{N} Peter Nyikos, \textsl{The Theory of Nonmetrizable Manifolds},
in K Kunen and J Vaughan, eds, ``Handbook of Set-Theoretic Topology" (Elsevier,
1984), 634--684.
\bibitem{Nyk} Peter Nyikos, \textsl{Various Smoothings of the Long Line and
Their Tangent Bundles}, Advances in Mathematics,
93(1992), 129--213.
\bibitem{Ny} Peter Nyikos, \textsl{Mary Ellen Rudin's Contributions to the
Theory of Nonmetrizable Manifolds}, in F D Tall, ed,
``The Work of Mary Ellen Rudin" (Annals of the New York Academy of
Sciences, 705(1993)), 92--113.
\bibitem{P} A R Pears, \textsl{Dimension Theory of General Spaces},
Cambridge University Press, 1975.
\bibitem{Po1} H Poincar\'e, \textsl{Second complement a l'analysis situs},
Proceedings of the London Mathematical Society, 32(1900), 277--302.
\bibitem{Po2} H Poincar\'e, \textsl{Cinquieme complement a l'analysis situs},
Rendiconti Circulo Matematici Palermo, 18(1904), 45--110.
\bibitem{ReZ} G M Reed and P L Zenor, \textsl{A metrization theorem for
normal Moore spaces}, in N M Stavrakas and K R Allen,
eds, ``Studies in Topology" (Academic Press, New York, 1974), 485--488.
\bibitem{ReZ2} G M Reed and P L Zenor, \textsl{Metrization of Moore spaces and
generalized manifolds}, Fundamenta Mathematicae, 91(1976), 203--210.
\bibitem{R} M E Rudin, \textsl{The undecidability of the existence of a
perfectly normal non-metrizable manifold}, Houston Journal of Mathematics,
5(1979), 249--252.
\bibitem{Ru} M E Rudin, \textsl{Two Nonmetrizable Manifolds}, Topology and
its Applications, 35(1990), 137--152.
\bibitem{RZ} M E Rudin and P L Zenor, \textsl{A perfectly normal
non-metrizable manifold}, Houston Journal of Mathematics, 2(1976), 129--134.
\bibitem{S} Yu M Smirnov, \textsl{On Metrization of Topological Spaces},
Uspehi Matematiki Nauk, 6(1951), \#6(46), 100--111 (American Mathematical Society
Translations \#91).
\bibitem{T} M G Tka\v{c}enko, \textsl{Ob odnom svoistve bicompactov (On a
property of compact spaces)}, Seminar po obshchei topologii (A seminar on general
topology), Moscow State University P. H., Moscow (1981), 149--156 (Russian).
\bibitem{WW} J M Worrell Jr and H H Wicke, \textsl{Characterizations of
Developable Topological Spaces}, Canadian Journal of Mathematics, 17(1965),
820--830.
\bibitem{WZ} Scott W Williams and Haoxuan Zhou, \textsl{Strong Versions of
Normality}, General Topology and its Applications, Proceedings of the 5th NE
Conference, New York 1989, in Lecture Notes in Pure and Applied Mathematics, 134
(Marcel Dekker, New York, 1991), 379--389.
\end{thebibliography}
\end{document}